\def\ver{Dec. 18, 2002, v.2}
\documentstyle{amsppt}
\magnification=1200
\hoffset=-10pt
\voffset=-30pt
\hsize=6.5truein
\vsize=8.9truein
\topmatter
\title Filtrations on Chow Groups\\
and Transcendence Degree
\endtitle
\author Morihiko Saito
\endauthor
\affil RIMS Kyoto University, Kyoto 606-8502 Japan \endaffil
\abstract For a smooth complex projective variety
$ X $ defined over a number field, we have filtrations
on the Chow groups depending on the choice of realizations.
If the realization consists of mixed Hodge structure
without any additional structure, we can show that the
obtained filtration coincides with the filtration of Green and
Griffiths, assuming the Hodge conjecture.
In the case the realizations contain Hodge structure and etale
cohomology, we prove that if the second graded piece of the
filtration does not vanish, it contains a nonzero element
which is represented by a cycle defined over a field of
transcendence degree one.
This may be viewed as a refinement of results of Nori, Schoen,
and Green-Griffiths-Paranjape.
For higher graded pieces we have a similar assertion assuming
a conjecture of Beilinson and Grothendieck's generalized Hodge
conjecture.
\endabstract
\endtopmatter
\tolerance=1000
\baselineskip=12pt
\def\mtim{\hbox{$\times$}}
\def\motm{\hbox{$\otimes$}}
\def\bC{{\Bbb C}}
\def\bP{{\Bbb P}}
\def\bQ{{\Bbb Q}}
\def\bZ{{\Bbb Z}}
\def\cD{{\Cal D}}
\def\cH{{\Cal H}}
\def\cM{{\Cal M}}
\def\cO{{\Cal O}}
\def\cP{{\Cal P}}
\def\oC{\overline{C}}
\def\ok{\overline{k}}
\def\oK{\overline{K}}
\def\oS{\overline{S}}

\def\tL{\widetilde{L}}
\def\tZ{\widetilde{Z}}
\def\CH{\hbox{{\rm CH}}}
\def\rank{\text{{\rm rank}}\,}
\def\supp{\text{{\rm supp}}\,}
\def\Spec{\text{{\rm Spec}}\,}
\def\Im{\hbox{{\rm Im}}}
\def\Ker{\hbox{{\rm Ker}}}
\def\Hom{\hbox{{\rm Hom}}}
\def\cHom{{{\Cal H}om}}
\def\Ext{\hbox{{\rm Ext}}}
\def\MHS{\text{{\rm MHS}}}
\def\Gal{\text{{\rm Gal}}}
\def\Gr{\hbox{{\rm Gr}}}
\def\alg{\text{\rm alg}}
\def\Mur{\text{\rm Mur}}

\def\SameName{\vrule height3pt depth-2.5pt width1cm}

\document
\centerline{{\bf Introduction}}

\bigskip\noindent
Let
$ X_{\bC} $ be a smooth complex projective variety, and
$ \CH^{p}(X_{\bC})_{\bQ} $ be the Chow group with rational
coefficients.
Choosing a category
$ \cM $ of realizations (see [12], [13], [20]), we can define a
filtration
$ F_{\cM} $ on
$ \CH^{p}(X_{\bC})_{\bQ} $ by spreading cycles out, see (1.4)
below (and also
$ [1], [16], [25], [29]) $.
By definition
$ F_{\cM}^{1}\CH^{p}(X_{\bC})_{\bQ} $ consists of null
homologous cycles, and
$ F_{\cM}^{2}\CH^{p}(X_{\bC})_{\bQ} $ is contained in the
kernel of the Abel-Jacobi map (tensored with
$ \bQ $).
It is conjectured that the filtration
$ F_{\cM} $ does not depend on the choice of
$ \cM $, and coincides with Murre's conjectural filtration [24].
We can verify this conjecture, assuming a conjecture of
Beilinson on the injectivity of the Abel-Jacobi map for smooth
projective varieties over number fields [2] together with the
Hodge conjecture.

In this paper we assume that
$ X_{\bC} $ is defined over a number field
$ k $.
Then a similar filtration has been defined by
M.~Green and P.~Griffiths [16], and we have

\medskip\noindent
{\bf 0.1.~Proposition.}
{\it If
$ \cM $ is the category of mixed Hodge structure without any
additional structure, then the filtration
$ F_{\cM} $ coincides with the filtration
$ F_{G} $ of Green and Griffiths {\rm [16],} assuming the Hodge
conjecture.
}

\medskip
Let
$ X $ be a smooth projective
$ k $-variety whose base change by
$ k \to \bC $ is
$ X_{\bC} $.
Let
$ K $ be a subfield of
$ \bC $ containing
$ k $,
and having finite transcendence degree.
Let
$ X_{K} $ be the base change of
$ X $ by
$ k \to K $.
Then
$ \CH^{p}(X_{K})_{\bQ} $ is identified with a subgroup of
$ \CH^{p}(X_{\bC})_{\bQ} $,
and has the induced filtration
$ F_{\cM} $.
It has been observed by Green and Griffiths [16] that the
property of this induced filtration is very much influenced by
the transcendence degree
$ d $ of
$ K $.
For example,
$ \Gr_{F_{\cM}}^{r}\CH^{p}(X_{K})_{\bQ} $ vanishes for
$ r > d + 1 $ if the realization consists of mixed Hodge
structure.
If
$ d = 0 $,
it is conjectured that
$ F_{\cM}^{2}\CH^{p}(X_{K})_{\bQ} = 0 $ by the above
conjecture of Beilinson.
However, for
$ d = 1 $,
it is shown by M.~Nori and C.~Schoen [31] that the kernel of
the Albanese map for certain surfaces has a nontrivial cycle
defined over a subfield
$ K $ of transcendence degree
$ 1 $.
Here we can show also the nonvanishing of
$ \Gr_{F_{\cM}}^{2}\CH^{2}(X_{\bC})_{\bQ} $ (see [25]),
which implies that the above estimate is optimal.
The results of Nori and Schoen are recently generalized by
Green-Griffiths-Paranjape [17] to the case of surfaces having
a nontrivial global
$ 2 $-form.
Considering these, we may have

\medskip\noindent
{\bf 0.2.~Conjecture.}
{\it If
$ \Gr_{F_{\cM}}^{r}\CH^{p}(X_{\bC})_{\bQ} \ne 0 $ with
$ r \ge 1 $, then it contains a nonzero element which is
represented by a cycle defined over a subfield of transcendence
degree
$ r - 1 $.
}

\medskip
In this paper we prove

\medskip\noindent
{\bf 0.3.~Theorem.}
{\it Assume that the realizations contain mixed Hodge structure
and \'etale cohomology with Galois action.
Then Conjecture {\rm (0.2)} is true for
$ r = 1, 2 $.
Assume further that Grothendieck's generalized Hodge conjecture
holds, and the filtration
$ F_{\cM} $ coincides with the filtration
$ F_{\MHS} $ associated to the category of realization
consisting of mixed Hodge structure.
Then Conjecture {\rm (0.2)} is true also for
$ r \ge 3 $.
}

\medskip
The proof uses Terasoma's argument on Hilbert's irreducibility
theorem [33] as in [17].
The same argument was also indicated by A.~Tamagawa when we
tried to construct an
$ l $-adic theory of normal functions [28].
It is quite interesting that we cannot prove Theorem (0.3) by
using only Hodge theory.
The hypothesis of (0.3) for
$ r = 2 $ is satisfied for
$ 0 $-cycles if
$ X $ has a nontrivial global
$ 2 $-from [29] (this follows from Murre's Albanese motive [23]
and Bloch's diagonal cycle [7]).
So Theorem (0.3) may be viewed as a refinement of the result
of Green-Griffiths-Paranjape [17].
For cycles of arbitrary codimension, we have a similar assertion
if the standard conjecture of Lefschetz-type holds for
$ X $.

If we restrict to the subgroup
$ \CH_{\alg}^{p}(X_{\bC})_{\bQ} $ consisting of cycles
algebraically equivalent to zero,
$ F_{\cM}^{2} $ coincides with the kernel of the Abel-Jacobi
map (or that of the
$ l $-adic Abel-Jacobi map).
This applies to the case of
$ 0 $-cycles on surfaces, and we have in general
$ \Gr_{F_{\cM}}^{r}\CH^{p}(X_{\bC})_{\bQ} = 0 $ for
$ r > p $ (see also [16]).
However it is not yet clear whether the nonvanishing of the
kernel of the Albanese map for a surface
$ X_{\bC} $ implies that
$ \Gr_{F_{\cM}}^{2}\CH^{2}(X_{\bC})_{\bQ} \ne 0 $, because
it is not proved that the filtration
$ F_{\cM} $ is separated.

I would like to thank L.~Barbieri-Viale, A.~Rosenschon and
A.~Tamagawa for useful discussions.

\bigskip\bigskip\centerline{{\bf 1. Filtrations on Chow groups}}

\bigskip\noindent
{\bf 1.1.~Realizations.}
We will denote by
$ \cM $ a category of (systems of) realizations, see [12], [13],
[20], etc.
The simplest example in our case is the abelian category
$ \MHS $ of
$ \bQ $-mixed Hodge structures whose graded pieces
$ \Gr_{m}^{W} $ are polarizable [11].
In this paper we choose a number field
$ k $ contained in
$ \bC $.
Then we have the category
$ \MHS_{k} $ of mixed
$ \bQ $-Hodge structures with
$ k $-structure, see [25], [29], etc.
Let
$ \ok $ be the algebraic closure of
$ k $ in
$ \bC $,
and put
$ G = \Gal(\ok/k) $.
For a prime number
$ l $,
we have an abelian category
$ \cM_{l} $ whose object consists of filtered vector spaces
$ (H_{\bQ},W) $ over
$ \bQ $,
$ (H_{l},W) $ over
$ \bQ_{l} $ and
$ (H_{\bC},F) $ over
$ \bC $ together with isomorphisms
$$
\alpha_{l} : (H_{\bQ},W)\motm_{\bQ}\bQ_{l}= (H_{l},W),\quad
\alpha_{\bC} : H_{\bQ}\motm_{\bQ}\bC = H_{\bC},
$$
where
$ (H_{l},W) $ is endowed with a continuous action of
$ G $ and
$ (\Gr_{m}^{W}H_{\bQ},\Gr_{m}^{W}(H_{\bC},F)) $ is a polarizable
$ \bQ $-Hodge structure of weight
$ m $ for any
$ m $ (here
$ W $ denotes also the induced filtration on
$ H_{\bC}) $,
see [12], [13], [20], etc.
We assume that polarizations are compatible with the Galois
action to assure the semisimplicity of pure objects.

For other examples, we have
$ \cM_{\text{\'et}} $ by considering
$ H_{l} $ for any prime numbers
$ l $,
and
$ \cM_{k,l}, \cM_{k,\text{\'et}} $ by considering also the
$ k $-structure.
It is also possible to consider the category of systems of
realizations as in [20].

Note that the category
$ \cM $ can be extended naturally to the category of mixed
sheaves
$ \cM(S/k) $ for any
$ k $-variety
$ S $,
and there is a forgetful functor from
$ \cM(S/k) $ to the category of perverse sheaves [6], see [29]
for the details.

\medskip\noindent
{\bf 1.2.~Deligne cohomology.}
Let
$ \cM $ be one of the categories of realizations as in (1.1).
Let
$ X $ be a smooth
$ k $-variety.
Then the cohomology
$ H^{i}(X/k,\bQ) $ is well-defined in
$ \cM $, using de Rham cohomology of
$ X_{k} $, \'etale cohomology of
$ X_{\ok} $, and cohomology of
$ X_{\bC} $ together with comparison isomorphisms,
see [12], [13], [20], etc.
Furthermore, there exists canonically
$ K_{\cH}(X/k) $ in the bounded derived category
$ D^{b}\cM $ whose cohomology is isomorphic to the cohomology of
$ X $ (using, for example, two sets of affine open coverings
associated to general hyperplane sections [5], see also
[25, 1.1]).

We define Deligne cohomology by
$$
H_{\cD}^{i}(X/k,\bQ(j)) = \Hom_{D^{b}\cM}(\bQ,K_{\cH}(X/k)(j)
[i]),
$$
where
$ (j) $ is the Tate twist, and
$ [i] $ is the shift of complexes.
If
$ \cM = \MHS $,
it is called the absolute
$ p $-Hodge cohomology in [3].
If
$ \cM = \MHS $ or
$ \MHS_{k} $,
then higher extension groups vanish in
$ \cM $ as a corollary of [10] (see [29]), and we have
canonical short exact sequences
$$
\aligned
0 \to \Ext_{\cM}^{1}(\bQ,H^{i-1}(X/k,\bQ)(j))
&\to H_{\cD}^{i}(X/k,\bQ(j))
\\
&\quad \to \Hom_{\cM}(\bQ,H^{i}(X/k,\bQ)(j)) \to 0.
\endaligned
\leqno(1.2.1)
$$
For a closed subvariety
$ Z $ of
$ X $,
we can define similarly the Deligne local cohomology
$ H_{\cD,Z}^{i}(X/k,\bQ(j)) $ using a complex
$ K_{\cH,Z}(X/k) $,
which is the shifted mapping cone of
$ K_{\cH}(X/k) \to K_{\cH}((X\setminus Z)/k) $.

We have the cycle map
$$
cl : \CH^{p}(X)_{\bQ} \to H_{\cD}^{2p}(X/k,\bQ(p)),
\leqno(1.2.2)
$$
which is compatible with the usual cycle class map to
$ H^{2p}(X_{\bC},\bQ)(p) $.
Its restriction to the null homologous cycles coincides with
Griffiths' Abel-Jacobi map [18] tensored with
$ \bQ $ if
$ k = \bC $,
$ \cM = \MHS $ and
$ X $ is smooth proper, see [9], [14], [15], [19], etc.
We can show that (1.2.2) is compatible with the direct image by
a proper morphism and the pull-back by any morphism, and
hence with the action of a correspondence, cf. [27].

\medskip\noindent
{\bf 1.3.~Leray filtration.}
Let
$ X, S $ be a smooth
$ k $-varieties.
Then the Deligne cohomology
$ H_{\cD}^{i}(X\mtim_{k}S/k,\bQ(j)) $ has the (decreasing)
Leray filtration
$ F_{L} $ induced by the canonical filtration
$ \tau $ on
$ K_{\cH}(X/k) $ using the canonical isomorphism
$$
K_{\cH}(X\mtim_{k}S/k) = K_{\cH}(X/k)\motm K_{\cH}(S/k).
$$
Here
$ F_{L}^{r} $ on
$ H_{\cD}^{i}(X\mtim_{k}S/k,\bQ(j)) $ is induced by
$ \tau_{\le i-r} $ as in [11].
Assume
$ X $ is smooth proper.
Then the filtration
$ F_{L} $ splits because we have a non canonical isomorphism
$$
K_{\cH}(X/k) \simeq \sum_{j} H^{j}(X/k,\bQ)[-j]\quad
\text{in}\,\, D^{b}\cM.
\leqno(1.3.1)
$$
(This follows from a general property of pure complexes,
see e.g. [26].)
In particular, for a morphism
$ S' \to S $,
the filtration
$ F_{L} $ is strictly compatible with the pull-back morphism
$$
H_{\cD}^{i}(X\mtim_{k}S/k,\bQ(j)) \to
H_{\cD}^{i}(X\mtim_{k}S'/k,\bQ(j)).
$$

By the canonical filtration on
$ K_{\cH}(S/k) $,
we have for each
$ m \in \bZ $ the Leray spectral sequence
$$
\aligned
E_{2}^{p,q}
&= \Ext_{\cM}^{p-m}(\bQ,H^{m}(X/k,\bQ)\motm
H^{q}(S/k,\bQ)(j))
\\
&\Rightarrow \Gr_{F_{L}}^{p+q-m}
H_{\cD}^{p+q}(X\mtim_{k}S/k,\bQ(j))
\endaligned
\leqno(1.3.2)
$$
It is conjectured that this degenerates at
$ E_{2} $, because
$ K_{\cH}(S/k) $ would be defined in the (conjectural) category
of motives where higher extension groups should vanish so that
a decomposition similar to (1.3.1) would hold.

We will denote by
$ F'_{L} $ the decreasing filtration on
$ \Gr_{F_{L}}^{r}H_{\cD}^{i}(X\mtim_{k}S/k,\bQ(j)) $ induced
by the canonical filtration
$ \tau $ on
$ K_{\cH}(S/k) $ so that
$ \Gr_{F'_{L}}^{s}\Gr_{F_{L}}^{r}H_{\cD}^{i}(X\mtim_{k}S/k,
\bQ(j)) $ is a subquotient of
$ \Ext_{\cM}^{s}(\bQ,H^{i-r}(X/k,\bQ)\motm H^{r-s}(S/k,
\bQ)(j)) $.

In the case
$ \cM = \MHS $ or
$ \MHS_{k} $,
the higher extension groups really vanish so that (1.3.2)
degenerates at
$ E_{2} $ and we get canonical
short exact sequences
$$
\aligned
0 \to \Ext_{\cM}^{1}(\bQ,H^{i-r}(X/k,\bQ)\motm H^{r-1}(S/k,
\bQ)(j)) \to \Gr_{F_{L}}^{r}H_{\cD}^{i}(X\mtim_{k}S/k,\bQ(j))
&
\\
\to \Hom_{\cM}(\bQ,H^{i-r}(X/k,\bQ)\motm H^{r}(S/k,\bQ)(j))
\to 0.
&
\endaligned
\leqno(1.3.3)
$$
In particular,
$ F_{L}^{\prime 2}\Gr_{F_{L}}^{r} = 0 $ in this case.

\medskip\noindent
{\bf 1.4.~Filtration on Chow groups.}
Let
$ X $ a smooth
$ k $-variety, and
$ K $ be a subfield of
$ \bC $ containing
$ k $,
and having finite transcendence degree over
$ k $.
Put
$ X_{K} = X\motm_{k}K $,
and
$ X_{\bC} = X\motm_{k}\bC $.
Then we have natural injections
$$
\CH^{p}(X_{K})_{\bQ} \to \CH^{p}(X_{\bC})_{\bQ},
\leqno(1.4.1)
$$
and
$ \cup_{K} \CH^{p}(X_{K})_{\bQ} = \CH^{p}(X_{\bC})_{\bQ} $.

Let
$ \zeta \in \CH^{p}(X_{K})_{\bQ} $.
By spreading out [7], there exists an irreducible smooth affine
$ k $-variety
$ S $ such that
$ k(S) = K $ and
$ \zeta $ is defined over
$ S $,
i.e.
there exists
$ \zeta_{S} \in \CH^{p}(X\mtim_{k}S)_{\bQ} $ whose
restriction to
$ X_{K} $ is
$ \zeta $,
where
$ X_{K} $ is identified with the generic fiber of
$ X\mtim_{k}S $
$ \to S $.
For an open subvariety
$ S' $ of
$ S $,
let
$ \zeta_{S'} $ denote the restriction of
$ \zeta_{S} $ over
$ S' $.
Then the limit of
$ \zeta_{S'} $ is well-defined, see [7].

Let
$ k_{S} $ be the algebraic closure of
$ k $ in
$ \Gamma (S,\cO_{S}) $,
and put
$ S_{\bC} = S\motm_{k_{S}}\bC $.
This is an irreducible variety, i.e.
$ S $ is geometrically irreducible over
$ k_{S} $.
(If we consider
$ S\motm_{k}\bC $ instead of
$ S\motm_{k_{S}}\bC $,
then the former is a disjoint union of copies of the latter
in the case
$ k_{S} $ is a normal extension of
$ k $.)
Note that
$ X\mtim_{k}S = X_{k_{S}}\mtim_{k_{S}}S $, and this allows
us to replace
$ k $ with
$ k_{S} $.
Actually we can replace
$ k $ with any finite extension, because we take the limit over
$ K $.

The cycle map (1.2.2) induces
$$
cl : \CH^{p}(X\mtim_{k}S)_{\bQ} \to H_{\cD}^{2p}
(X\mtim_{k}S/k_{S},\bQ(p)),
\leqno(1.4.2)
$$
and the filtration
$ F_{\cM} $ on
$ \CH^{p}(X\mtim_{k}S)_{\bQ} $ is defined to be the
induced filtration by the Leray filtration
$ F_{L} $ on
$ H_{\cD}^{2p}(X\mtim_{k}S/k_{S},\bQ(p)) $.
Then, taking the inductive limit over the non empty open
subvarieties of
$ S $,
we get the filtration
$ F_{\cM} $ on
$ \CH^{p}(X_{K})_{\bQ} $.

This means that
$ \zeta \in F_{\cM}^{r}\CH^{p}(X_{K})_{\bQ} $ if
$ cl(\zeta_{S}) \in F_{L}^{r}H_{\cD}^{2p}(X\mtim_{k}S/k_{S},
\bQ(p)) $ for some
$ S $,
and hence
$ \Gr_{F_{\cM}}^{r}\zeta $ is nonzero in
$ \Gr_{F_{\cM}}^{r}\CH^{p}(X_{K})_{\bQ} $ if the restrictions of
$ \Gr_{F_{L}}^{r}cl(\zeta_{S}) $ to
$ \Gr_{F_{L}}^{r}H_{\cD}^{2p}(X\mtim_{k}S'/k_{S},\bQ(p)) $ does
not vanish for any non empty open subvarieties
$ S' $ of
$ S $.

We can show that
$ F_{\cM} $ is strictly compatible with the base change by
$ K \to K' $,
see [29].
This implies that
$ \CH^{p}(X)_{\bQ} $ has the filtration
$ F_{\cM} $ which is strictly compatible with (1.4.1).

\medskip\noindent
{\bf 1.5.~Filtration of Green and Griffiths.}
In the case
$ \cM = \MHS $, a similar filtration is constructed by
M.~Green and P.~Griffiths [16].
They assume that the
$ S $ in (1.4) are smooth {\it projective}, and then,
roughly speaking, consider everything modulo ambiguity
coming from cycles over proper closed subvarieties of
$ S $ (here they also assume Grothendieck's generalized Hodge
conjecture).
More precisely, for a smooth projective
$ k $-variety
$ S $ and a divisor
$ Z $ of
$ S $ defined over
$ k_{S} $,
we have an exact sequence
$$
\CH^{p-1}(X\mtim_{k}Z)_{\bQ} \to \CH^{p}(X\mtim_{k}
S)_{\bQ} \to \CH^{p}(X\mtim_{k}(S\setminus Z))_{\bQ}
\to 0,
\leqno(1.5.1)
$$
and the filtration
$ F_{G} $ of Green and Griffiths on
$ \CH^{p}(X\mtim_{k}(S\setminus Z))_{\bQ} $ is defined
to be the quotient filtration of
$ F_{\cM} $ on
$ \CH^{p}(X\mtim_{k}S)_{\bQ} $,
where
$ \cM = \MHS $.
Then we take the inductive limit as before.

\medskip\noindent
{\bf 1.6.~Proposition.}
{\it
$ F_{G} = F_{\cM} $, assuming the Hodge conjecture.
}

\medskip\noindent
{\it Proof.}
It is enough to show the assertion on
$ \CH^{p}(X\mtim_{k}(S\setminus Z))_{\bQ} $.
This is reduced to the case
$ Z $ is a divisor with normal crossings by using an embedded
resolution.
We have a canonical morphism of (1.5.1) to
$$
H_{\cD,Z}^{2p}(X\mtim_{k}S/k_{S},\bQ(p)) \to H_{\cD}^{2p}
(X\mtim_{k}S/k_{S},\bQ(p)) \to H_{\cD}^{2p}(X\mtim_{k}
(S\setminus Z)/k_{S},\bQ(p)).
$$
Here we may assume
$ k_{S} = k $ (and similarly for intersections of
irreducible components of
$ Z $) replacing
$ k $ if necessary.
Assuming the Hodge conjecture, we have to prove the following:

\medskip
For
$ \zeta \in F_{\cM}^{r}\CH^{p}(X\mtim_{k}S)_{\bQ} $
such that
$ \Gr_{F_{L}}^{r}cl(\zeta) \in \Gr_{F_{L}}^{r}
H_{\cD}^{2p}(X\mtim_{k}S,\bQ(p)) $ comes from
$ \xi \in \Gr_{F_{L}}^{r}H_{\cD,Z}^{2p}(X\mtim_{k}S,
\bQ(p)) $,
there exists
$ \zeta' \in \CH^{p-1}(X\mtim_{k}Z)_{\bQ} $ such that
the image of
$ cl(\zeta') $ in
$ H_{\cD}^{2p}(X\mtim_{k}S,\bQ(p)) $ belongs to
$ F_{L}^{r} $,
and coincides with
$ \Gr_{F_{L}}^{r}cl(\zeta) $ modulo
$ F_{L}^{r+1} $.

\medskip
This is verified by using correspondences
$ \Gamma_{a} \in \CH^{\dim S-1}(S\mtim_{k}\tZ)_{\bQ} $ such
that
$$
(\Gamma_{a})_{*} : H^{j}(S,\bQ) \to H^{j-2}(\tZ,\bQ)(-1)
$$
vanishes for
$ j \ne a $,
and the restriction of
$ i_{*}(\Gamma_{a})_{*} $ to
$ \Im \,i_{*} \subset H^{a}(S,\bQ) $ is the identity for
$ j = a $,
where
$ \tZ $ is the normalization of
$ Z $.
Indeed, if we denote by
$$
\xi_{0} \in \Hom_{\MHS}(\bQ,H^{2p-r}(X,\bQ)\motm
H_{Z}^{r}(S,\bQ)(p))
$$
the image of
$ \xi $ by the canonical morphism, then the Hodge conjecture
implies the existence of
$ \zeta' \in \CH^{p-1}(X\mtim_{k}\tZ)_{\bQ} $
such that the K\"unneth component of the cycle class of
$ \zeta' $ in
$ \Hom(\bQ,H^{2p-a}(X,\bQ)\motm H_{Z}^{a}(S,\bQ)(p))
$ coincides with
$ \xi_{0} $ for
$ a = r $,
and is zero otherwise.
We may assume further that the image of
$ cl(\zeta') $ in
$ H_{\cD}^{2p}(X\mtim_{k}S,\bQ(p)) $ belongs to
$ F_{L}^{r} $ by modifying
$ \zeta' $ using
$ \Gamma_{a} $ for
$ a < r $ together with the decomposition (1.3.1).
So the assertion is reduced to the case
$ \xi_{0} = 0 $ by modifying
$ \zeta $ using
$ \zeta' $.
Then the assertion follows by using
$ \Gamma_{a} $ for
$ a = r $.

\bigskip\bigskip
\centerline{{\bf 2. Proof of Theorem (0.3)}}

\bigskip\noindent
{\bf 2.1.~Hilbert's irreducibility theorem.}
We first recall Terasoma's argument [33] on Hilbert's irreducibility
theorem, which is essential for the proof of (0.3).
Let
$ U $ be a non empty open subvariety of
$ \bP_{k}^{1} $,
and
$$
0 \to L \to \tL \to \bQ_{l,U} \to 0
\leqno(2.1.1)
$$
be a short exact sequence of smooth
$ \bQ_{l} $-sheaves on
$ U $,
where
$ \bQ_{l,U} $ denotes the constant sheaf of rank one on
$ U $.
Put
$ K = k(U) $, and let
$ \oK $ be an algebraic closure of
$ K $.
There exists a
$ k $-valued point
$ x $ of
$ U $ such that (2.1.1) splits if and only if its restriction
over
$ x $ does.
Indeed, choosing a geometric point over
$ x $ on each Galois \'etale covering of
$ U $ in a compatible way with natural projections, we get a
morphism of
$ \Gal(\ok/k) $ to
$ \pi_{1}(U,\Spec \oK) $,
and hence to the arithmetic monodromy group of
$ \tL $.
Then we have infinitely many
$ k $-valued points
$ x $ such that the last morphism is surjective by Hilbert's
irreducibility theorem [22] together with the structure of the
$ l $-adic monodromy group [32], see [33].
Related to the
$ l $-adic theory of normal functions, the same argument was
indicated by A.~Tamagawa, see [28].

Here it is also possible to get infinitely many
$ k $-valued points
$ x $ such that the above property holds for the monodromy
groups of
$ L $ and
$ \tL $ simultaneously by the theory of Hilbert set.
Note also that the exact sequence (2.1.1) can be replaced
by a short exact sequence
$ 0 \to L^{1} \to \tL \to L^{0} \to 0 $ of smooth
$ \bQ_{l} $-sheaves, because
$$
\Ext^{1}(L^{0},L^{1}) =
\Ext^{1}(\bQ_{l,U},\cHom(L^{0},L^{1})).
$$

\medskip\noindent
{\bf 2.2.~Restriction of extension classes.}
Let
$ f : S \to U $ be a smooth projective morphism of smooth
irreducible
$ k $-varieties where
$ U $ is a non empty open subvariety of
$ \bP_{k}^{1} $.
Let
$ n = \dim S - 1 $, and
$ L = R^{n}f_{*}\bQ_{X} \in \cM(U/k) $ where
$ \cM(U/k) $ denotes the category of mixed sheaves on
$ U $ (shifted by
$ \dim U) $, and
$ L $ is pure of weight
$ n $, see [29].
Here we assume that there is a forgetful functor from
$ \cM $ to
$ \cM_{l} $ in (1.1).
By semisimplicity we have a direct sum decomposition
$$
L = L' \oplus L''\quad \text{in}\,\,\cM(U/k)
$$
such that
$ H^{0}(U/k,L') = 0 $ and
$ L'' $ is constant over
$ \Spec k $ (i.e. the pull-back of an object on
$ \Spec k $ by the structure morphism).

Let
$ H $ be a pure object of weight
$ n + 1 $ in
$ \cM $ (e.g.
a direct factor of
$ H^{i}(X/k,\bQ)(q) $ for a smooth projective
$ k $-variety
$ X $ where
$ i - 2q = n + 1 $).
Let
$ H_{U} = {a}_{U}^{*}H $,
where
$ a_{U} : U \to \Spec k $ is the structure morphism.
By the adjunction for
$ a_{U} $,
we have a natural isomorphism
$$
\Ext_{\cM(U/k)}^{1}(H_{U},L) =
\Hom_{D^{b}\cM}(H,(a_{U})_{*}L[1]).
$$
This implies
$$
\aligned
\Ext_{\cM(U/k)}^{1}(H_{U},L')
&= \Hom_{\cM}(H,H^{1}(U/k,L')),
\\
\Ext_{\cM(U/k)}^{1}(H_{U},L'')
&= \Ext_{\cM}^{1}(H,L''_{x}),
\endaligned
\leqno(2.2.1)
$$
for any
$ k $-valued point
$ x $ of
$ U $,
because
$ \Hom_{\cM}(H,H^{1}(U/k,L'')) = 0 $.

Let
$ \xi \in \Hom_{\cM}(H,H^{1}(U/k,L')) $.
The corresponding extension class is denoted also by
$ \xi $.
If
$ \xi \ne 0 $,
there exists a
$ k $-valued point
$ x $ of
$ U $ such that the restriction
$ \xi_{x} $ of
$ \xi $ to
$ x $ does not vanish by (2.1), because (2.2.1) holds also for
$ l $-adic sheaves.
Note that the same argument still holds after replacing
$ k $ by a finite extension.
In the case
$ \dim S = 1 $ and
$ n = 0 $,
we may also assume that
$ f^{ -1}(x) $ consists of one point.
Then replacing
$ U, L $ with
$ S, \bQ_{l,S} $, the restriction of
$ \xi $ to some
$ k $-valued point of
$ S $ does not vanish (replacing
$ k $ if necessary).

\medskip\noindent
{\bf 2.2.~Restriction to open subvarieties.}
With the above notation and assumptions, let
$ S_{x} = f^{-1}(x) $.
Then
$ L'_{x} $ is a direct factor of
$ H^{n}(S_{x}/k,\bQ) $,
and we get
$$
\xi_{x} \in \Ext_{\cM}^{1}(H, H^{n}(S_{x}/k,\bQ)).
$$
We now consider to restrict
$ \xi_{x} $ to a non empty open subvariety
$ S'_{x} $ of
$ S_{x} $.
We assume that the underlying Hodge structure of
$ H $ does not have a nontrivial subobject with level
$ < n $, where the level of a Hodge structure is the
difference between the maximal and minimal numbers
$ p $ such that
$ \Gr_F^{p} \ne 0 $
(and the difference between level and weight is even).
Let
$$
H' = H^{n}(S'_{x}/k,\bQ).
$$
It has weights
$ \ge n $.
If
$ S'_{x} $ is sufficiently small,we have
$$
W_{n}H' = H^{n}(S_{x}/k,\bQ)/N^{1}H^{n}(S_{x}/k,\bQ),
$$
where
$ N $ is the `coniveau' filtration.
By semisimplicity there exists a subobject
$ H'' $ such that
$$
H^{n}(S_{x}/k,\bQ) = N^{1}H^{n}(S_{x}/k,\bQ) \oplus H''.
$$

We have
$ \Hom_{\cM}(H,H'/W_{n}H')
= \Hom_{\cM}(H,\Gr_{n+1}^{W}H') = 0 $, because
$ \Gr_{n+1}^{W}H' $ has level
$ < n $ (see [11]).
Then, using the long exact
sequence associated to
$$
0 \to W_{n}H' \to H' \to H'/W_{n}H' \to 0,
\leqno(2.3.1)
$$
we see that the restriction of
$ \xi_{x} $ to
$ S'_{x} $ does not vanish if
$ \xi_{x} $ does not come from
$ \Ext_{\cM}^{1}(H,N^{1}H^{n}(S_{x}/k,\bQ)) $
(i.e. if its image in
$ \Ext_{\cM}^{1}(H,H'') $ does not vanish).

In the case
$ \dim S = 2 $ and
$ n = 1 $, the last condition is trivially satisfied
because
$ N^{1}H^{1}(S_{x}/k,\bQ) = 0 $.
Furthermore,
$ H'/W_{n}H' $ is a direct sum of copies of
$ \bQ $, replacing
$ k $ with a finite extension (depending on
$ S'_{x}) $ if necessary.
Indeed, it is given by taking a basis of the kernel of the
cycle class map
$ \sum_{i} \bZ[D_{i}] \to H^{2}(S_{x}/k,\bQ)(1) $
where the
$ D_{i} $ are the irreducible components of
$ S_{x} \setminus S'_{x} $, which may be assumed to be
absolutely irreducible (replacing
$ k $ if necessary).
This fact will be used in (2.4).

In general,
$ L'_{1,x} := N^{1}H^{n}(S_{x}/k,\bQ) \cap L'_{x} $ does
not vanish.
However, it corresponds to a
$ \bQ_{l} $-submodule stable by the action of
$ \Gal(\ok/k) $, and is hence extended to an \'etale subsheaf
$ L'_{1} $ of
$ L' $ by the argument in (2.1).
Let
$ s = \rank L'_{1} $.
Then taking the pull-back to
$ U_{\bC} $, it determines a subsheaf with
$ \bQ $-coefficients, and the latter underlies a variation of
Hodge structure.
Indeed,
$ \wedge^{s}L'_{1} $ determines a variation of Hodge
structure of rank
$ 1 $ contained in
$ \wedge^{s}L' $ by the global invariant cycle theorem
(using a finite covering if necessary, because the monodromy of
$ \wedge^{s}L'_{1} $ is defined over
$ \bZ $ and is finite, see [11]).
Then
$ L'_{1} $ is the kernel of
$ L' \to \wedge^{s+1}L' $ defined locally by a generator of
$ \wedge^{s}L'_{1} $, and hence underlies a variation of
Hodge structure.

This argument implies that the restriction of
$ \xi_{x} $ to
$ S'_{x} $ does not vanish if
$ \xi \in \Hom_{\cM}(H,H^{1}(U/k,L')) $ is nonzero.
Indeed, if the restriction vanishes, the corresponding
$ l $-adic extension class comes from
$ L'_{1,x} \,(\subset L'_{x}) $ which is extended to
$ L'_{1} \,(\subset L') $.
We apply some argument in (2.1) also to
$ L'/L'_{1} $, where we may assume
$ H = \bQ $ by the last remark of (2.1) and the monodromy
group of the extension for
$ L'/L'_{1} $ is a quotient of that for
$ L' $.
Then we see that
$ \xi $ comes from
$ \Hom(H,H^{1}(U/k,L'_{1})) $, where
$ \Hom $ is considered in
$ \cM_{l} $.
But
$ H^{1}(U/k,L'_{1}) $ has level
$ < n $,
because the stalk of
$ L'_{1} $ has level
$ \le n - 2 $, see [34].
So
$ \xi $ vanishes by the hypothesis on the level of
$ H $, and the assertion follows.

\medskip\noindent
{\bf 2.4.~Proof of (0.3).}
By hypothesis there exists a smooth irreducible affine
$ k $-variety
$ S $ together with
$ \zeta \in \CH^{p}(X\mtim_{k}S)_{\bQ} $ such that its
cycle class
$ cl(\zeta) $ in
$ H_{\cD}^{2p}(X\mtim_{k}S/k,\bQ(p)) $ belongs to
$ F_{L}^{r} $,
and the restriction of
$ \Gr_{F_{L}}^{r}cl(\zeta) $ to
$ \Gr_{F_{L}}^{r}H_{\cD}^{2p}(X\mtim_{k}S'/k,\bQ(p)) $
does not vanish for any non empty open subvariety
$ S' $ of
$ S $.
Using the spectral sequence (1.3.2),
$ \Gr_{F_{L}}^{r}cl(\zeta) $ induces
$$
\xi_{0} \in \Hom_{\cM}(\bQ,H^{2p-r}(X/k,\bQ)\motm H^{r}
(S'/k,\bQ)(p)).
$$

We first consider the case where
$ \xi_{0} $ does not vanish for any
$ S' $.
Let
$ d = \dim X - p $.
Since
$ \xi_{0} $ corresponds to the morphism
$$
\xi'_{0} : H^{2d+r}(X/k,\bQ)(d) \to H^{r}(S'/k,\bQ)),
$$
this nonvanishing is equivalent to that the image of
$ \xi'_{0} $ has level
$ r $ (assuming Grothendieck's generalized Hodge conjecture for
$ r > 2 $).
So we may assume
$ \dim S = r $ by the weak Lefschetz theorem.
Put
$ n = r - 1 $.
Let
$$
H = H^{2d+r}(X/k,\bQ)(d),
$$
and
$ H_{<n} $ be the largest subobject of
$ H $ which has level
$ < n $.
By semisimplicity there exists a subobject
$ H_{> n} $ with a decomposition
$ H = H_{<n} \oplus H_{> n} $, and the restriction of
$ \xi'_{0} $ to
$ H_{> n} $ does not vanish.
So the assertion follows from (2.2-3) applied to a
Lefschetz pencil.

Now we may assume
$ \xi_{0} = 0 $, i.e.
$ \Gr_{F_{L}}^{r}cl(\zeta) \in F_{L}^{\prime 1}
\Gr_{F_{L}}^{r} $, see (1.3).
Then
$ \Gr_{F_{L}}^{r}cl(\zeta) $ induces
$$
\aligned
\xi_{1}
&\in \Ext_{\cM}^{1}(\bQ,H^{2p-r}(X/k,\bQ)\motm
H^{r-1}(S'/k,\bQ)(p))
\\
&= \Ext_{\cM}^{1}(H,H^{r-1}(S'/k,\bQ)).
\endaligned
$$
Consider the case where
$ \xi_{1} $ does not vanish for any
$ S' $.
If
$ r = 1 $, we may replace
$ S' $ with any point (replacing
$ k $ if necessary), and the assertion is clear.
So we may assume
$ r > 1 $.
In this case we have to show the nonvanishing of its
restriction to any non empty open subvariety
$ C' $ of a general hyperplane section
$ \oC $ of a smooth projective compactification
$ \oS $ of
$ S $.

If
$ r = 2 $, let
$ \cP_{\oS/k}, \cP_{\oC/k} $ be the
Picard variety of
$ \oS, \oC $.
Then we have an injective morphism of
$ k $-varieties
$ \cP_{\oS/k} \to \cP_{\oC/k} $,
and any
$ k $-valued point on the image can be lifted to a
$ k $-valued point of
$ \cP_{\oS/k} $.
So the assertion follows using the
short exact sequence (2.3.1) for
$ S' $ and
$ C' $ (and replacing
$ k $ if necessary).

If
$ r > 2 $,
we may assume Grothendieck's generalized Hodge conjecture, and
the `coniveau' filtration
$ N $ coincides with the filtration by the level of Hodge
structure.
If
$ S' $ is a sufficiently small open affine subvarieties of
$ \oS $,
then
$$
W_{n}H^{n}(S'/k,\bQ) = H^{n}(\oS/k,\bQ)/N^{1}
H^{n}(\oS/k,\bQ),
$$
and similarly for
$ C' $.
By the weak Lefschetz theorem, the restriction morphism
$$
H^{n}(\oS/k,\bQ) \to H^{n}(\oC/k,\bQ)
$$
is injective, and splits by semisimplicity.

Assume that the pull-back of
$ \xi_{1} $ by
$ C' \to S' $ vanishes.
Then, using the long exact sequence associated to (2.3.1), we
see that
$ \xi_{1} $ factors through a direct factor of
$ H $ (or equivalently, of
$ H^{2p-r}(X/k,\bQ)) $ with level
$ < n $,
because
$ \Gr_{n+1}^{W}H^{n}(C'/k,\bQ) $ has level
$ < n $.
By the Hodge conjecture there exists a smooth proper
$ k $-variety
$ Y $ of pure dimension
$ r - 2 $ together with a correspondence
$ \Gamma \in \CH^{p-1}(Y\mtim_{k}X)_{\bQ} $ such that
the image of
$$
\Gamma_{*} : H^{r-2}(Y/k,\bQ) \to H^{2p-r}(X/k,\bQ)(p-r+1)
$$
coincides with
$ N^{p-r+1}H^{2p-r}(X/k,\bQ)(p-r+1) $ (i.e.
the maximal subobject with level
$ \le r -2) $.
We have also a correspondence
$ \Gamma' \in \CH^{\dim X-p+r-1}(X\mtim_{k}Y)_{\bQ} $
such that the restriction of
$ \Gamma_{*}\Gamma'_{*} $ to
$ \Im \,\Gamma_{*} \subset H^{2p-r}(X/k,\bQ) $ is the
identity.
So we may replace
$ \zeta $ with
$ \Gamma_{*}\Gamma'_{*}\zeta $ to show the vanishing of
$ \xi_{1} $.
Here
$ \zeta $ is extended to
$ X\mtim_{k}\oS $ by taking the closure, and the
correspondences preserve the filtration
$ \tau $ because they induce morphisms of complexes
$ K_{\cH}(X/k) $, see [27].
Since
$ \Gamma' $ induces
$$
\Gamma'_{*} : \CH^{p}(X\mtim_{k}\oS) \to
\CH^{r-1}(Y\mtim_{k}\oS),
$$
we see that
$ \supp \Gamma'_{*}\zeta \subset Y\mtim_{k}Z $ with
$ Z $ a divisor on
$ \oS $,
because
$ r - 1 > \dim Y $.
So we get the assertion, because
$ \supp \Gamma_{*}\Gamma'_{*}\zeta \subset X\mtim_{k}Z $.

Now we may assume further
$ \xi_{1} = 0 $, i.e.
$ \Gr_{F_{L}}^{r}cl(\zeta) \in F_{L}^{\prime 2}
\Gr_{F_{L}}^{r} $.
If
$ r > 2 $, we have
$ \Gr_{F_{L}}^{r}cl(\zeta) = 0 $ by the hypothesis on the
coincidence of the two filtrations, because
$ F_{L}^{\prime 2}\Gr_{F_{L}}^{r} = 0 $ for
$ \cM = \MHS $ and the filtrations
$ F_{\cM} $ and
$ F'_{\cM} $ in (1.4) are functorial for
$ \cM $.
So we may assume
$ r = 2 $, since the case
$ r = 1 $ is trivial by the vanishing of
$ H^{r-2}(S'/k,\bQ) $.
Then
$ \Gr_{F_{L}}^{r}cl(\zeta) $ induces
$$
\xi_{2} \in \Ext_{\cM}^{2}(\bQ,H^{2p-2}(X/k,\bQ)\motm
H^{0}(S'/k,\bQ)(p)),
$$
because
$ d_{2} : E_{2}^{m,1} \to E_{2}^{m+2,0} $ vanishes in (1.3.2)
(replacing
$ k $ if necessary) where
$ m = 2p-2, j = p $.
Indeed,
$ H^{0}(S'/k,\bQ) = \bQ $ is a direct factor of
$ K_{\cH}(S'/k) $ by choosing a
$ k $-valued point
$ x $ of
$ S' $,
because we have
$ \bQ \to K_{\cH}(S'/k) \to \bQ $ by the structure morphism
and
$ x $.
In this case, the assertion is clear because
$ H^{0}(S'/k,\bQ) = H^{0}(C'/k,\bQ) = \bQ $ (replacing
$ k $ if necessary).
Thus we have verified all the cases, because
$ \Gr_{F_{L}}^{r}cl(\zeta) = 0 $ if
$ \xi_{2} = 0 $ and
$ r = 2 $.
This completes the proof of Theorem (0.3).

\medskip\noindent
{\bf 2.5.~Remarks.} (i)
It is conjectured that the filtration
$ F_{\cM} $ is separated, and gives the conjectural
``motivic'' filtration of Beilinson [4] and Bloch [7].
This depends on the injectivity of the Abel-Jacobi map for
smooth projective
$ k $-varieties, which is also a conjecture of Beilinson [2],
see also [8], [16], [29], [30], etc.
It is expected that the filtration
$ F_{\cM} $ does not depend on the choice of
$ \cM $,
and coincides with Murre's (conjectural) filtration
$ F_{\Mur} $ [24].
Indeed, we have
$$
F_{\Mur} \subset F_{\cM}\quad \text{and}\quad
F_{\Mur} = F_{\cM} \mod \cap _{i} F_{\cM}^{i},
\leqno(2.5.1)
$$
see [29].
The existence of
$ F_{\Mur} $ can be deduced from the separatedness of
the filtration
$ F_{\cM} $ assuming the algebraicity of the K\"unneth
components of the diagonal, see [21].
The separatedness of
$ F_{\cM} $ is reduced to the above conjecture of
Beilinson on the Abel-Jacobi map for
$ k $-varieties, assuming the Hodge conjecture in the
case the codimension of cycles is more than
$ 2 $.

\medskip
(ii) We have
$ \Gr_{F_{\cM}}^{r}\CH^{p}(X_{\bC})_{\bQ} = 0 $ for
$ r > p $.
If
$ \cM = \MHS $ or
$ \MHS_{k} $,
then
$ \Gr_{F_{\cM}}^{r}\CH^{p}(X_{K})_{\bQ} = 0 $ for
$ r > \text{\rm tr deg}\, K/k + 1 $.
These follow from the vanishing of
$ H^{i}(S_{\bC},\bQ) $ for a smooth affine variety
$ S $ and
$ i > \dim S $, together with the compatibility of the
cycle map with the pull-back by a closed embedding
(and the vanishing of higher extension groups).
These assertions have been shown by M.~Green and P.~Griffiths
[16] for their filtration, assuming the above conjecture of
Beilinson and Grothendieck's generalized Hodge conjecture.
Note that these conjectures imply also that the filtration is
separated and ends at the
$ p $-th step.

\medskip
(iii) Restricted to the subgroup
$ \CH_{\alg}^{p}(X)_{\bQ} $ consisting of cycles
algebraically equivalent to
$ 0 $, the kernel of the Abel-Jacobi map coincides with
$ F_{\cM}^{2}\CH_{\alg}^{p}(X_{\bC})_{\bQ} $,
see [30], 3.9.
Indeed, for a curve
$ C $ and a correspondence
$ \Gamma \in \CH^{p}(C\mtim X)_{\bQ} $,
we have a decomposition
$ H^{1}(C,\bQ) = \Im \,\Gamma_{*} \oplus \Ker \,\Gamma_{*} $
induced by idempotents of
$ \CH^{1}(C\mtim C)_{\bQ} $,
where
$ \Gamma_{*} : H^{1}(C,\bQ) \to H^{2p-1}(X,\bQ)(p-1) $.
By a similar argument, the kernel of the usual Abel-Jacobi map
coincides with that of the
$ l $-adic Abel-Jacobi map on
$ \CH_{\alg}^{p}(X)_{\bQ} $ (because we have the
injectivity in the divisor case using the Kummer sequence).

\bigskip\bigskip
\centerline{{\bf References}}

\bigskip

\item{[1]}
M.~Asakura, Motives and algebraic de Rham cohomology,
in: The arithmetic and geometry of algebraic cycles (Banff),
CRM Proc. Lect. Notes, 24, AMS, 2000, pp. 133--154.

\item{[2]}
A.~Beilinson, Higher regulators and values of
$ L $-functions, J. Soviet Math. 30 (1985), 2036--2070.

\item{[3]}
A.~Beilinson, Notes on absolute Hodge cohomology, Contemporary
Math. 55 (1986) 35--68.

\item{[4]}
\SameName, Height pairing between algebraic cycles, Lect. Notes
in Math., vol. 1289, Springer, Berlin, 1987, pp. 1--26.

\item{[5]}
\SameName, On the derived category of perverse sheaves, ibid.
pp. 27--41.

\item{[6]}
A.~Beilinson, J.~Bernstein and P.~Deligne, Faisceaux pervers,
Ast\'erisque, vol. 100, Soc. Math. France, Paris, 1982.

\item{[7]}
S.~Bloch, Lectures on algebraic cycles, Duke University
Mathematical series 4, Durham, 1980.

\item{[8]}
\SameName, Algebraic cycles and values of
$ L $-functions, J. Reine Angew. Math. 350 (1984), 94--108.

\item{[9]}
\SameName, Algebraic cycles and the Beilinson conjectures,
Contemporary Math. 58 (1) (1986), 65--79.

\item{[10]}
J.~Carlson, Extensions of mixed Hodge structures, in
Journ\'ees de G\'eom\'etrie Alg\'ebrique d'Angers 1979,
Sijthoff-Noordhoff Alphen a/d Rijn, 1980, pp. 107--128.

\item{[11]}
P.~Deligne, Th\'eorie de Hodge I, Actes Congr\`es Intern.
Math., 1970, vol. 1, 425-430; II, Publ. Math. IHES, 40 (1971),
5--57; III ibid., 44 (1974), 5--77.

\item{[12]}
\SameName, Valeurs de fonctions L et p\'eriodes d'int\'egrales,
in Proc. Symp. in pure Math., 33 (1979) part 2, pp. 313--346.

\item{[13]}
P.~Deligne, J.~Milne, A.~Ogus and K.~Shih, Hodge Cycles,
Motives, and Shimura varieties, Lect. Notes in Math., vol 900,
Springer, Berlin, 1982.

\item{[14]}
F.~El Zein and S.~Zucker,
Extendability of normal functions associated to algebraic
cycles, in Topics in transcendental algebraic geometry, Ann.
Math. Stud., 106, Princeton Univ. Press, Princeton, N.J., 1984,
pp. 269--288.

\item{[15]}
H.~Esnault and E.~Viehweg, Deligne-Beilinson cohomology, in
Beilinson's conjectures on Special Values of
$ L $-functions, Academic Press, Boston, 1988, pp. 43--92.

\item{[16]}
M.~Green and P.~Griffiths, Hodge-theoretic invariants for
algebraic cycles, Intern. Math. Res. Notices 2003 (9), 477--510.

\item{[17]}
M.~Green, P.~Griffiths and K.H.~Paranjape, Cycles over fields
of transcendence degree one, preprint (math.AG/0211270).

\item{[18]}
P.~Griffiths, On the period of certain rational integrals I,
II, Ann. Math. 90 (1969), 460--541.

\item{[19]}
U.~Jannsen, Deligne homology, Hodge-$ D $-conjecture, and
motives, in Beilinson's conjectures on Special Values of
$ L $-functions, Academic Press, Boston, 1988, pp. 305--372.

\item{[20]}
\SameName, Mixed motives and algebraic
$ K $-theory, Lect. Notes in Math., vol. 1400, Springer, Berlin,
1990.

\item{[21]}
\SameName, Motivic sheaves and filtrations on Chow groups,
Proc. Symp. Pure Math. 55 (1994), Part 1, pp. 245--302.

\item{[22]}
S.~Lang, Fundamentals of Diophantine geometry, Springer, Berlin,
1983.

\item{[23]}
J.~Murre, On the motive of an algebraic surface, J. Reine
Angew. Math. 409 (1990), 190--204.

\item{[24]}
\SameName, On a conjectural filtration on Chow groups of an
algebraic variety, Indag. Math. 4 (1993), 177--201.

\item{[25]}
A.~Rosenschon and M.~Saito, Cycle map for strictly decomposable
cycles, preprint (math.AG/0104079).

\item{[26]}
M.~Saito, Mixed Hodge Modules, Publ. RIMS, Kyoto Univ., 26 (1990),
221--333.

\item{[27]}
\SameName, Hodge conjecture and mixed motives, I, Proc. Symp.
Pure Math. 53 (1991), 283--303; II, in Lect. Notes in Math.,
vol. 1479, Springer, Berlin, 1991, pp. 196--215.

\item{[28]}
\SameName, Tate conjecture and mixed perverse sheaves, preprint.

\item{[29]}
\SameName, Arithmetic mixed sheaves, Inv. Math. 144 (2001),
533--569.

\item{[30]}
\SameName, Refined cycle maps, in Alg. Geom. 2000 (Azumino),
Adv. Stud. Pure Math 36, 2002, pp. 115--143.

\item{[31]}
C.~Schoen, Zero cycles modulo rational equivalence for some
varieties over fields of transcendence degree one, Proc. Symp.
Pure Math. 46 (1987), part 2, pp. 463--473.

\item{[32]}
J.P.~Serre, Groupes analytiques
$ p $-adiques, S\'eminaire Bourbaki, 270 (1963/64).

\item{[33]}
T.~Terasoma, Complete intersections with middle Picard number
$ 1 $ defined over {\bf Q}, Math. Z. 189 (1985), 289--296.

\item{[34]}
S. Zucker, Hodge theory with degenerating coefficients,
$ L_{2} $-cohomology in the Poincar\'e metric,
Ann. Math., 109 (1979), 415--476.

\medskip\noindent
\ver

\bye